# On Ethnomathematics: In Memory of Ubiratan D'Ambrosio

Veselin Jungic
Department of Mathematics
Simon Fraser University

The purpose of this note is to honour memory of Ubiratan D'Ambrosio, a Brazilian mathematics educator, who passed away on May 12, 2021.

I also wish to encourage the reader to read (or re-read) D'Ambrosio's article "Ethnomathematics and its place in the history and pedagogy of mathematics," in *For the Learning of Mathematics* (D'Ambrosio, 1985).

This, by far the most frequently cited article among D'Ambrosio's 107 publications listed on the Sematic Scholar website[1], has left a long-lasting mark on how the mathematics community understands the connection between mathematics and underlying cultures.

It should be mentioned that over the years, D'Ambrosio and his collaborators further developed the ideas and concepts related to ethnomathematics introduced in (D'Ambrosio, 1985). See, for example, (D'Ambrosio, 2006; D'Ambrosio, 2017; Rosa & Orey, 2021; Scott, 2012; Shannon, 2021) and the references within.

This note is inspired by two models presented in (D'Ambrosio, 1985, pp. 45-46) that were based on the idea that ethnomathematics, or any other knowledge or way of knowing, is a human way of interacting with reality as an individual and as a society. I will present brief versions of those two models.

## Ethnomathematics

D'Ambrosio defined ethnomathematics as "the mathematics which is practiced among identifiable cultural groups, such as national tribal societies, labour groups, children of certain age bracket, professional classes, and so on" (D'Ambrosio, 1985, p. 45). In contrast, he defined academic mathematics as "mathematics which is taught and learned in schools." In addition, D'Ambrosio argued that mathematics includes "apart from the Platonic ciphering and arithmetic, mensuration and relation of planetary orbits, the capabilities of classifying, ordering, inferring and modelling."

In short, D'Ambrosio's concept of ethnomathematics includes:

> […] a very broad range of human activities which, throughout history, have been expropriated by the scholarly establishment, formalized and codified and incorporated into

---
[1] https://www.semanticscholar.org/author/Ubiratan-D%E2%80%99ambrosio/116668173

what we call academic mathematics. But which remain alive in culturally identified groups and constitute routines in their practices. (D'Ambrosio, 1985, p. 45)

As an example that supports D'Ambrosio's claim that ethnomathematics "remains alive," i.e., something that is active and alert, I highlight the work of a contemporary Coast Salish artist Qwul'thilum (DylanThomas) of the Lyackson First Nation.

In (Thomas & Schattschneider, 2011), Qwul'thilum describes his 2007 print entitled *Sacred Cycle*, a remarkable piece of Salish ethnomathematics:

> Sacred Cycle depicts three identical salmon in a circle, each related to the other by a $120^0$ rotation. The eyes of the salmon contain simple faces, a design element I learned from Rande[2]. A salmon repeats three times to complete the cycle then it starts over again. Unlike Escher, who would have carved a woodblock with just one fish and then rotated the block to produce a print with rotation symmetry, my silkscreen for this print had to contain the complete 3-fold cycle of salmon. (Thomas & Schattschneider, 2011, p. 202)

In the same article, Qwul'thilum and his coauthor Doris Schattschneider, an American mathematician, explain that Qwul'thilum's early work was partly inspired by the artist Susan Point of the Musqueam First Nation:

> The design of *Sacred Cycle* was heavily influenced by the works of Susan Point, who often uses rotational symmetry in her designs. It is due to her great influence on the Salish art community that geometric symmetry has become common in contemporary Salish art. (Thomas & Schattschneider, 2011, p. 202)

The fact that the Salish art community accepted *geometric symmetry* as a common component of their work connects with one of the aims of what D'Ambrosio called *the Project Ethnomathematics*:

> understanding how, historically, societies absorb innovation, how do they incorporate the new into their quotidian. (D'Ambrosio, 2006, 82)

Qwul'thilum's work also illustrates how ethnomathematics, as a manifestation of a certain culture, "is subject to inter and intra-cultural encounters" (D'Ambrosio, 2006, p. 76).

For example, Qwul'thilum's print from 2008 entitled *Salmon Spirits* (Thomas & Schattschneider, 2011, p. 204, Figure 5) was inspired by tessellations created by M.C. Escher, a Dutch graphic artist:

> After 6 months of apprenticeship, I felt that my skills in traditional Salish design were strong enough to undertake an Escher-style tessellation using Salish design. The problem was that I had no further education in mathematics past grade 12 and had no idea how to create a tessellation. I spent days staring at Escher's symmetry drawings such as his

---

[2] Rande Cook, a Kwakwaka'wakw multimedia artist, was Qwul'thilum's mentor at the time.

butterfly and lizard tessellations and from these, soon figured out how I could make a tessellation. (Thomas & Schattschneider, 2011, p. 202)

Schattschneider, who has extensively written about the work of M.C. Escher (see, for example, (Schattschneider, 2017) and the references within) observes:

It is interesting to note that [Qwul'thilum's] path to making his first tessellation parallels that of Escher. Escher, too, had no mathematical background for this task, and figured it out by studying several geometric tessellations by majolica tiles in the Alhambra. (Thomas & Schattschneider, 2011, p. 203)

Hence, Qwul'thilum's print *Salmon Spirits*, undoubtedly a piece of Salish art, and thus Salish ethnomathematics, is, in part, an outcome of a sequence of inter-cultural encounters drawing back to medieval artifacts of Islamic culture. But *inter-cultural encounters* also include those when one cultural group imposes its set of value on another group:

We should not forget that colonialism grew together in a symbiotic relationship with modern science, in particular with mathematics, and technology. (D'Ambrosio, 1985, p. 47)

## A Forest of Ethnomathematics

In 1990, *For the Learning of Mathematics* published an article entitled "Ethnomathematics and Education" by Marcelo C. Borba, another Brazilian mathematics educator. (Borba, 1990)

Borba offers the following description of ethnomathematics:

Mathematical knowledge expressed in the language code of a given sociocultural group is called "ethnomathematics." (Borba, 1990, p. 40)

Furthermore, Borba reflects that the ethno- in ethnomathematics "should be understood as referring to cultural groups, and not to the anachronistic concept of race," and that mathematics "should be seen as a set of activities such as ciphering, measuring, classifying, ordering, inferring and modelling." (Borba, 1990, p. 40)

Seemingly diverging from D'Ambrosio's concept of academic mathematics, Borba observes that "academic mathematics is not universal (in the sense of being independent of culture)" and states:

even the mathematics produced by professional mathematicians can be seen as a form of ethnomathematics because it was produced by an identifiable cultural group and because it is not the only mathematics that has been produced. (Borba, 1990, p. 40)

This is something that resonates well with the experience of the author of this note, a professional mathematician. Consequently, I paraphrase Borba to obtain the following statement:

Even knowledge that belongs to a specific area of mathematics, like combinatorics, analysis, experimental mathematics, or mathematical education, for example, can be seen as a form of ethnomathematics because it is produced by an identifiable cultural group and because it is not the only mathematics that has been produced.

Here I use the term *cultural group* in the sense of D'Ambrosio's definition:

> A group of individuals is identified as a cultural group if it reveals shared knowledge and compatible behavior, both subordinated to a set of values. (D'Ambrosio, 2006, 77)

As a piece of evidence that academics involved in a specific area of mathematics form a cultural group, in the sense of D'Ambrosio's definition, I offer the following example.

Attendees of mathematical meetings that organize sessions from different areas, like the semiannual Canadian Mathematical Society meetings, for example, likely join a session from their own area of interest. There they share their knowledge with others interested in a similar set of mathematical problems or mathematics related issues. They socialize with their past, present, or future collaborators, students, teachers, or supervisors, i.e., people with compatible academic and non-academic behavior. They are among colleagues who submit their work to, act as referees for, or are editors of the academic journals that are setting the values for the area.

Borba's suggestion that a "forest might be a better image [than a straight line] of the whole set of ethnomathematics, in which each tree would be considered as a different expression of ethnomathematics, socio-culturally produced" (Borba, 1990, p. 41) inspired the visualization of ethnomathematics depicted in Figure 1.

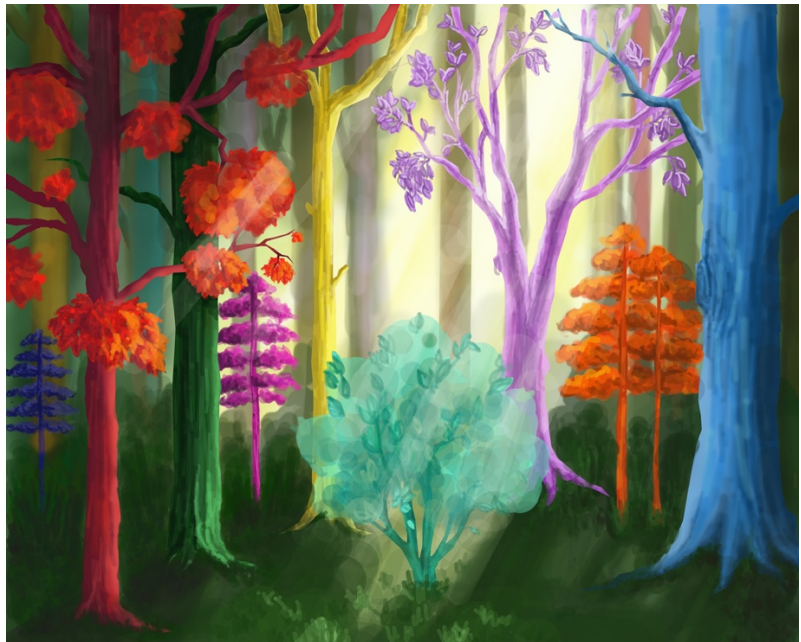

*Figure 1    Forest of ethnomathematics by Veselin Jungic and Nicole Oishi*

## Where Does Ethnomathematics Come From?

D'Ambrosio (1985, pp. 45-46) describes the "ceaseless cycle"

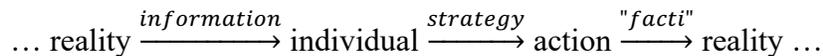
$$\ldots \text{reality} \xrightarrow{information} \text{individual} \xrightarrow{strategy} \text{action} \xrightarrow{"facti"} \text{reality} \ldots$$

as "the basis for the theoretical framework upon which we base our ethnomathematics concept."

Next, D'Ambrosio observes that "individual behavior is homogenized in certain ways […] to build up societal behavior" and obtains the following model

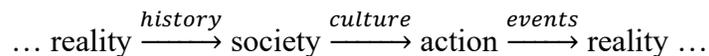
$$\ldots \text{reality} \xrightarrow{history} \text{society} \xrightarrow{culture} \text{action} \xrightarrow{events} \text{reality} \ldots$$

Furthermore, D'Ambrosio explains:

> […] culture manifests itself through jargons, codes, myths, symbols, utopias, and ways of reasoning and inferring. Associated with these we have practises such as ciphering and counting, measuring, classifying, ordering, inferring, modelling, and so on, which constitute ethnomathematics. (D'Ambrosio, 1985, p. 46)

The following quote by Borba is a restatement of the main message that emerges from D'Ambrosio's models, i.e., the view that ethnomathematics is a human way of interacting with reality.

> The ethnomathematics of a cultural group is part of the group's life; the mathematics is generated by the culture in an "umbilical" way. Ethnomathematics is developed by the cultural group's interest in its problematic situations, which then further develop the group's interest in its ethnomathematics. This interest in ethnomathematics is natural because it was generated by the member of the cultural group in response to their own situations. (Borba, 1990, p. 41)

By including the idea that every tree in the forest of ethnomathematics is rooted in reality, we obtain the visualization of ethnomathematics that is presented in Figure 2[3].

---

[3] The reader will probably notice that the root of each tree in the forest of ethnomathematics is a connected graph with no cycles, hence a graph theoretical tree itself.

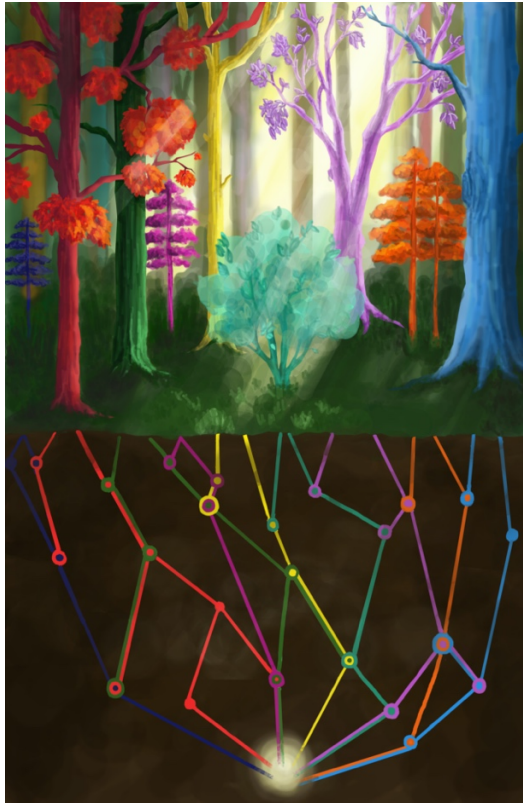

*Figure 2  Forest of ethnomathematics and its root by Veselin Jungic and Nicole Oishi*

## Conclusion

Ubiratan D'Ambrosio introduced the notion of ethnomathematics and tirelessly worked on its promotion and the further development over many years. This included acknowledging the dark side of mathematics, particularly its role in the process of colonization.

Still, D'Ambrosio remained positive and optimistic:

> Ethnomathematics is, thus, associated with the pursuit of PEACE. The main goal of Ethnomathematics is to build up a civilization free of truculence, arrogance, intolerance, discrimination, inequity, bigotry, and hatred. (D'Ambrosio, 2017, p. 665)

## References:

Borba, M. (1990) Ethnomathematics and education. *For the Learning of Mathematics* **10**(1), 39-43.

D'Ambrosio, U. (1985) Ethnomathematics and its place in the history and pedagogy of mathematics. *For the Learning of Mathematics* **5**(1), 44-48.

D'Ambrosio, U. (2006) The Program of Ethnomathematics and the challenges of globalization. *Circumscribere: International Journal for the History of Science* **1**, 74-82.

D'Ambrosio, U. (2017) Ethnomathematics and the pursuit of peace and social justice. *Educação Temática Digital* **19**(3), 653-666.

Rosa, M. & Orey, D. C. (2021) Ubiratan D'Ambrosio: Celebrating His Life and Legacy. *Journal of Humanistic Mathematics* **11**(2), 430-450. Available at:
https://scholarship.claremont.edu/jhm/vol11/iss2/26

Schattschneider, D. (2017) Lessons in Duality and Symmetry from M.C. Escher. In Fenyvesi, K. & Lähdesmäki, T, (Eds.) Aesthetics of Interdisciplinarity: Art and Mathematics, pp. 105–118. Cham, Switzerland: Springer.

Scott, P. (2012) The Intellectual Contributions of Ubiratan D'Ambrosio to Ethnomathematics. *Cuadernos de Investigación y Formación en Educación Matemática* **7**(10) 241-246

Shannon. A.G. (2021) Ubiratan D'Ambrosio [1932-2021] – ethnomathematics educator for the twenty-first century. *International Journal of Mathematical Education in Science and Technology* **52**(8), 1139-1142. https://doi.org/10.1080/0020739X.2021.1948629

Thomas, D. & Schattschneider, D. (2011) Dylan Thomas: Coast Salish artist. *Journal of Mathematics and the Arts* **5**(4), 199-211. https://doi.org/10.1080/17513472.2011.625346